# K-NOMIAL COEFFICIENTS

Florentin Smarandache
University of New Mexico
200 College Road
Gallup, NM 87301, USA
E-mail: smarand@unm.edu

In this article we will widen the concepts of ”binomial coefficients” and ”trinomial coefficients” to the concept of ”k-nomial coefficients”, and one obtains some general properties of these. As an application, we will generalize the” triangle of Pascal”.

Let’s consider a natural number $k \geq 2$; let $P(x) = 1 + x + x^2 + ... + x^{k-1}$ be the polynomial formed of k monomials of this type; we’ll call it ”k-nomial”.

We will call *k-nomial coefficients* the coefficients of the power of $x$ of $(1 + x + x^2 + ... + x^{k-1})^n$, for $n$ positive integer. We will note them $Ck_n^h$ with $h \in \{0,1,2,...,2pn\}$.

In continuation one will build by recurrence a triangle of numbers which will be called ” *triangle of the numbers of order k*”.

CASE 1: $k = 2p + 1$.

On the first line of the triangle one writes 1 and one calls it ”line 0”.

(1) It is agreed that all the cases which are to the left and to the right of the first (respectively of the last) number of each line will be consider like being 0. The lines which follow are called ”line 1”, ”line 2”, etc… Each line will contain $2p$ numbers to the left of the first number, $p$ numbers on the right of the last number of the preceding line. Numbers of the line $i+1$ are obtained by using those of the line $i$ in the following way:

$Ck_{i+1}^j$ is equal to the addition of $p$ numbers which are to its left on the line $i$ and of $p$ numbers which are to the right on the line $i$, to the number which is above it (see. Fig. 1). One will take into account the convention 1.

**Fig. 1**

$$
\begin{array}{llll}
 & \overbrace{\ldots\ldots\ldots}^{p\ \text{numbers}} & \cdot & \overbrace{\ldots\ldots\ldots}^{p\ \text{numbers}} \\
\text{line } i & & & \\
\text{line } i+1 & & \cdot Ck_{i+1}^{j} &
\end{array}
$$

Example for $k = 5$:

1

1 1 1 1 1

1 2 3 4 5 4 3 2 1

1 3 6 10 15 18 19 18 15 10 6 3 1

1 4 10 20 35 52 68 80 85 80 68 52 35 20 10 4 1

………………..………………………………………………………………

The number

$$C5_1^0 = 0 + 0 + 0 + 0 + 1 = 1;$$
$$C5_1^3 = 0 + 1 + 0 + 0 + 0 = 1;$$
$$C5_2^3 = 0 + 1 + 1 + 1 + 1 = 4;$$
$$C5_3^7 = 4 + 5 + 4 + 3 + 2 = 18;$$

etc.

**Properties of the triangle of numbers of order k:**

1) The line $i$ has $2pi + 1$ elements.
2) $Ck_n^h = \sum_{i=0}^{2p} Ck_{n-1}^{h-i}$ where by convention $Ck_n^t = 0$ for $\begin{cases} t < 0 \\ t > 2pr \end{cases}$ and

   This is obvious taking into account the construction of the triangle.
3) Each line is symmetrical relative to the central element.
4) First elements of the line $i$ are 1 and $i$.
5) The line $i$ of the triangle of numbers of order $k$ represent the k-nomial coefficients of $(1 + x + x^2 + ... + x^{k-1})^i$.

   The demonstration is done by recurrence on $i$ of $\mathbb{N}^*$:

   a) For $i = 1$ it is obvious; (in fact the property would be still true for $i = 0$).

   b) Let's suppose the property true for $n$. Then

$$(1 + x + x^2 + ... + x^{k-1})^{n+1} = (1 + x + x^2 + ... + x^{k-1})(1 + x + x^2 + ... + x^{k-1})^n =$$
$$= (1 + x + x^2 + ... + x^{2p}) \cdot \sum_{j=0}^{2pn} Ck_n^j \cdot x^j =$$
$$= \sum_{t=0}^{2p(n+1)} \sum_{\substack{i+j=t \\ 0 \le j \le 2p \\ 0 \le i \le 2pn}} Ck_n^i \cdot x^i \cdot x^j =$$
$$= \sum_{t=0}^{2p(n+1)} \left( \sum_{j=0}^{2p} Ck_n^{t-j} \right) x^t = \sum_{t=0}^{2p(n+1)} Ck_{n+1}^t \cdot x^t .$$

6) The sum of the elements locate on line $n$ is equal to $k^n$.

The first method of demonstration uses the reasoning by recurrence. For $n=1$ the assertion is obvious. One supposes the property truth for $n$, i.e. the sum of the elements located on the line $n$ is equal to $k^n$. The line $n+1$ is calculated using the elements of the line $n$. Each element of the line $n$ uses the sum which calculates each of $p$ elements locate to its left on the line $n+1$, each of $p$ elements locate to its right on the line $n+1$ and that which is located below: thus it is used to calculate $k$ numbers of the line $n+1$.

Thus the sum of the elements of the line $n+1$ is $k$ times larger than the sum of those of the line $n$, therefore it is equal to $k^{n+1}$.

7) The difference between the sum of the k-nomial coefficients of an even rank and the sum of the k-nomial coefficients of an odd rank located on the same line $\left(Ck_n^0 - Ck_n^1 + Ck_n^2 - Ck_n^3 + ...\right)$ is equal to 1.

One obtains it if in $(1+x+x^2+...+x^{k-1})^n$ one takes $x=-1$.

8) $Ck_n^0 \cdot Ck_m^h + Ck_n^1 \cdot Ck_m^{h-1} + ... + Ck_n^h \cdot Ck_m^0 = Ck_{n+m}^h$

This results from the fact that, in the identity

$$(1+x+x^2+...+x^{k-1})^n \cdot (1+x+x^2+...+x^{k-1})^m = (1+x+x^2+...+x^{k-1})^{n+m}$$

the coefficient of $x^h$ in the member from the left is $\sum_{i=0}^{h} Ck_n^i \cdot Ck_m^{h-i}$ and that of $x^h$ on the right is $Ck_{n+m}^h$.

9) The sum of the squares of the k-nomial coefficients locate on the line $n$ is equal to the k-nomial coefficient located in the middle of the line $2n$.

For the proof one takes $n=m=h$ in the property 8. One can find many properties and applications of these k-nomial coefficients because they widen the binomial coefficients whose applications are known.

CASE 2: $k=2p$.

The construction of the triangle of numbers of order $k$ is similar:

On the first line one writes 1; it is called line 0

The lines which follow are called line 1, line 2, etc. Each line will have $2p-1$ elements more than the preceding one; because $2p-1$ is an odd number, the elements of each line will be placed between the elements of the preceding line (which is different from the case 1 where they are placed below).

The elements locate on the line $i+1$ are obtained by using those of the line $i$ in the following way:

$Ck_{i+1}^j$ is equal to the sum of $p$ elements located to its left on the line $i$ with $p$ elements located to its right on the line $i$.

Fig. 2

line $i$ $\qquad \overbrace{\ldots\ldots\ldots}^{p\ numbers} \quad \cdot \quad \overbrace{\ldots\ldots\ldots}^{p\ numbers}$

line $i+1$ $\qquad\qquad \cdot Ck_{i+1}^{j}$

Example for $k = 4$:

1
1 1 1 1
1 2 3 4 3 2 1
1 3 6 10 12 12 10 6 3 1
1 4 10 20 31 40 44 40 31 20 10 4 1
……………………………………………………………….

From the property 1': $Ck_n^h = \sum_{i=0}^{2p-1} Ck_{n-1}^{h-i}$

By joining together properties 1 and 1': $Ck_n^h = \sum_{i=0}^{k-1} Ck_{n-1}^{h-i}$

The other properties of Case 1 are preserved in Case 2, with similar profs. However in the property 7, one sees that the difference between the sum of the k-nomial coefficients of even rank and that of the k-nomial coefficients of odd rank locate on the same line is equal to 0.